\documentclass{amsart}
\usepackage{amsfonts,amsmath,amssymb,eucal,tipa}
\usepackage[all]{xy}
\usepackage{hyperref}
\usepackage[utf8]{inputenc}

\begin{document}

\parindent=0pt
\parskip=6pt

\newcommand{\ie}{\textit{i}.\textit{e}.\,}
\newcommand{\eg}{\textit{e}.\textit{g}.\,}
\newcommand{\cf}{{\textit{cf}. \,}}
\newcommand{\F}{{\mathbb F}}
\newcommand{\C}{{\mathbb C}}
\newcommand{\bP}{{\mathbb P}}
\newcommand{\Q}{{\mathbb Q}}
\newcommand{\R}{{\mathbb R}}
\newcommand{\T}{{\mathbb T}}
\newcommand{\Z}{{\mathbb Z}}
\newcommand{\ch}{{\rm{ch}}}
\newcommand{\kt}{{\hat{K}}}
\newcommand{\pt}{{\rm pt}}
\newcommand{\vk}{{\varkappa}}
\newcommand{\tT}{{t_{\mathbb T}}}
\newcommand{\Sw}{{\mathcal S}}
\newcommand{\Spec}{{\rm Spec }}
\newcommand{\Spf}{{\rm Spf}}
\newcommand{\w}{{\mathfrak w}}
\newcommand{\li}{{\rm li}}
\newcommand{\half}{{\textstyle{\frac{1}{2}}}}
\newcommand{\ve}{{\varepsilon}}
\newcommand{\f}{{\mathbb F}}
\newcommand{\Td}{{\rm Td}}
\newcommand{\Wu}{{\rm Wu}}
\newcommand{\bro}{{\boldsymbol{\rho}}}
\newcommand{\btau}{{\boldsymbol{\tau}}}
\newcommand{\bmu}{{\boldsymbol{\lambda}}}
\newcommand{\PV}{{\sf pv}}
\newcommand{\spin}{{{\rm Spin}}}
\newcommand{\pl}{{ł}}

\newcommand{\bW}{{\overline{\mathbb W}}}
\newcommand{\bF}{{\overline{{\mathbb F}}_p}}
\newcommand{\bx}{{\bar{x}}}
\newcommand{\bw}{{\bar{w}}}
\newcommand{\bD}{{\hat{\mathbb D}}}
\newcommand{\sli}{{\widetilde{\rm{li}}}}
\newcommand{\hl}{{\widehat{\mathbb H}}}
\newcommand{\Zt}{{\Z_{(2)}}}

\newcommand{\dl}{{\ł}} 

\title{Periods for topological circle actions}

\author[Jack Morava]{Jack Morava}

\address{Department of Mathematics, The Johns Hopkins University,
Baltimore, Maryland}

\begin{abstract}{The language of currents and sparks \cite{5} may be useful in the study of Hopkin-Singer Wu classes \cite{9,22} in geometric topology, via distributional generalizations \cite{4,17} of classical zeta functions.}\end{abstract}

\maketitle \bigskip

{\bf \S I $\T$-equivariant Atiyah-Segal-Tate cohomology}\bigskip

{\bf 1.1} Atiyah-Segal equivariant $K$-theory $K^*_\T$ for the circle $\T \subset \C^\times$ takes values in $\Z_2$-graded modules over the representation ring
\[
K^*_\T(\pt) = R_\C(\T) = \Z[\vk,\vk^{-1}] \;,
\]
where $\vk = [\C(-1)]$ is the class of the one-dimensional complex representation of $\T$ in which $z \in \T$ acts as multiplication by $z^{-1}$. With this convention\begin{footnote}{opposite to that of \cite{1}(\S 2.6.6)}\end{footnote}, the Chern character 
\[
\ch : K^*_\T(\pt) \to H^{**}(B\T,\Q) \cong \Q[[c]]
\]
sends $\vk$ to $\exp(c)$, where $c \: (\; = [\frac{K(\nabla)}{2\pi i}]$) is the Euler class of the line bundle 
\[
E\T \times_\T \C(-1) \to E\T \times_\T \pt = B\T \cong \C P^\infty
\]
(measured by a Hermitian connection $\nabla$, \cite{16}(14.4, 14.10, Appendix C (Theorem)). Borel completion
\[
K^*_\T(X) \to K^*(X \wedge_\T E\T) = \kt^*_\T(X)
\]
sends 
\[
K^*_\T(\pt) = \Z[\vk^{\pm 1}] \to \kt^*_\T(\pt) = \Z[[q]] \;,
\]
where $q = 1 - \vk^{-1}$ with $\ch q = 1 - \exp(-c)$. The graded version of these constructions replaces $\Z$ with $K = \Z[v^{\pm 1}]$, where $v$ is a Bott element of cohomological degree -2, so $\vk \mapsto (1 - qv)^{-1} = \sum_{k \geq 0}(qv)^k$ . \bigskip

{\bf 1.2} The equivariant Tate cohomology $\tT^* \kt$ of the Borel completion fits in a long exact sequence of $\tT \kt$ - modules, which for a point becomes
\[
\xymatrix{
0 \ar[r] & \kt^*_\T (\pt) = \Z[[q]] \ar[r] & \tT^*\kt(\pt) \ar[r]^-\partial & K_*(B\T) \ar[r] & 0 \;,}
\]
where \cite{24}
\[
K_*(B\T) \cong \Z[b_k  \:|\: k \geq 1]/(b(z_0 + z_1 - z_0 z_1) = b(z_0)b(z_1))
\]
with $b(z) = 1 + \sum_{k \geq 1} b_k z^k = (1 - z)^b$, \ie $b_k \mapsto (-)^k \binom{b}{k} \in \Q[b]$. This defines a ring isomorphism $\tT^*\kt \cong q^{-1}\Z[[q]] = \Z((q))$ by setting
\[
\partial q^{-k} = \sum_{0 \leq j \leq k} Q^k_j \binom{b}{j}
\]
with $Q^k_j \in \Z, \; \eg Q^k_k = (-1)^k k!$. \bigskip

{\bf Proposition} {\it The $\T$-equivariant Tate cohomology $\tT^*K$ for Atiyah-Segal $K$-theory fits in a long exact sequence
\[
\xymatrix{
\dots \ar[r] & K^*_\T \ar[r] & \tT^*K \cong (1-\vk)^{-1} K^*_\T  \ar[r]^-\partial & K_*(B\T) \ar[r] & \dots}
\]
compatible with Borel completion, with $\tT^*K(\pt) = \Z[\vk^{\pm 1},(1- \vk)^{-1}]$.} 

The resulting cohomology theory can thus be interpreted \cite{18} as taking values in sheaves  of $\Z_2$ - graded modules on the projective line over the integers
\[
{\rm Spec} {} \tT K \cong \mathbb{P}^1(\Z) - \{0,1,\infty\}
\]
punctured at the three points defined over $\F_2$. We can ask about the behavior of this functor under automorphisms of the projective line.\bigskip

{\bf 1.3} In the diagram below, localizations and completions are as indicated by the notation, making free use of identities such as 
\[
(1-x)^{-1} = \sum_{k \geq 0} x^k, \; x^{-1} = \sum_{k \geq 0} (1-x)^k \; \dots 
\]
where available:
\[
\xymatrix{
{} & \Z[x] \ar[dl] \ar[d] \ar[ddr] \\
\Z[x,x^{-1}] \ar[d] & \Z[x, x(1-x)^{-1}] \ar[d] \\
\Z((x)) \ar[d] & \ar[l] \Z[x^{\pm 1}, (1-x)^{-1}] \ar[dr] \ar[dl] \ar[r] & \Z((1-x)) \ar[d] \\
\Q((c)) \ar@{<->}[rr]^-{\li_1}  & {} & \Q((h)) }
\]
Here $x \mapsto 1 - e^{-c}$ on the left, while $x \mapsto e^{-h}$ on the right.  

The dual diagram, in the language of spectral spaces, expresses the three-punctured sphere as something like the multiplicative group formally completed at its unit. The two formal points 
\[
\xymatrix{
\Spf \Q((c)) \ar[r] & \Spec \Z[x^{\pm 1}, (1-x)^{-1}]  & \ar[l] \Spf \Q((h)) }
\]
are witnessed by 
\[ 
1 - \exp(-c) \mapsto \exp(-h) ,  \ie  h = - \log |1 - \exp(-c)| :=  \li_1(-c) \;;
\]
but while 
\[
\rm{li}_1(x) = \log \circ \left[\begin{array}{cc}
                                                    1 & -1 \\
                                                     1 & 0
                                                     \end{array}\right] \circ \exp (x) \equiv - \log |x| 
\]
(modulo smooth functions of $x$) is real analytic on $\R_{>0}$ and defines a tempered distribution on $\R$ \cite{5}(\S 3 Cor 2), it is not representable by a formal power series. Thus $c$ and $h$ are not algebraically commensurable, but $h \to 0$ as $c \to \infty$ and $c \to 0$ as $h \to \infty$. In the theory of geometric quantization, this is an instance of the principle that waves at high frequency concentrate along geodesics. \bigskip

This suggests the interest of some kind of torsor \cite{10} 
\[
\Q - {\rm Alg} \ni A \mapsto {\rm NatTrans}_{\otimes A}(K_{\rm Tate},K_{\rm Borel}),                       
\]
of `transcendental' periods: monoidal equivalences between these two completed equivariant cohomology theories, regarded as a scheme of some kind over $\Q$.\bigskip

{\bf \S II Entropy, temperature, and action} \bigskip

{\bf 2.1 Wu's lemma}
\cite{20}(\S 2.1.3) asserts that 
\[
\half (e^{2x} - 1) \in \Z_{(2)}[[x]]
\] 
{\it has two-adic integral coefficients, and is congruent to $\sum_{k \geq 0} x^{2^k}$ modulo two.}\bigskip

Let $\bW$ be the ring of Witt vectors of $\bF$, and let $\bD$ be the completion $\bW \langle \langle F \rangle \rangle$ of the twisted power series algebra. Reduction modulo $p$ defines a ring homomorphism 
\[
\bW[[w]] \to \bF[[\bw]] \;;
\]
if we regard this as a morphism of  $\bD$-algebras by setting $F\bw = \bw^p$ with
\[
w(x) = p^{-1}(\exp(px) - 1) = x + \dots \in \Z_p[[x]], \; x = p^{-1} \log (1 + pw) \in \Z_p[[w]] \;;
\]
then
\[
\bw = \sum_{k \geq 0} \bx^{p^k} \;.
\]
defines ring isomorphisms $\Z_p[[w]] \cong \Z_p[[x]]$ and $\bF[\bw]] \cong \bF[[\bx]]$.
\bigskip

{\bf Proposition} $w \mapsto w - x := Fw = \half px^2 + \dots \in \Z_p[[w]] \;,$ 
\[ 
Fx = p^{-1} \log (\exp(px) - px) \in \Z_p[[x]]
\]
{\it defines a lift of Frobenius from $\bF[[\bw]]$ to $\bW[[w]]$ when} $p = 2$.  $\Box$ \bigskip

This seems to be more general than $x^p + p \delta(x)$, with $\delta$ a $p$-derivation in the sense of \cite{2}, see also \cite{3, 8}.\bigskip

{\bf Definition} {\it Wu's unit} is the ($\Zt$-adic) power series
\[
\w(x) := \frac{-2x}{1 - e^{-2x}} \in (\Z_{(2)}[[x]])^\times  \; \equiv \; (\sum_{k \geq 0} x^{2^k-1})^{-1} \in (\F_2[[x]])^\times \;.
\]
It is a deformation $\w \equiv x + \w^2$ of the identity object $\w(0) = 1$ in the pro\'etale stack of yes or no questions\begin{footnote}{\cf Hamming space $\prod^{\mathbb N}\F_2$
}\end{footnote}; in particular, $\w'(x) \equiv 1$ mod two.\bigskip

Let $L(x) = (1 + e^{-x})^{-1}$ be the  classical logistic ($0 \to 1)$ function, and in the language of Charles Epstein \cite{4}(Prop 1)) and myself, let us write
\[
\li_1(x) := - \log |1 - e^x|, \; \li_0(x) = - (1 - e^{-x})^{-1}, 
\]
with 
\[
d\li_1(x) = \li_0(x) \cdot dx = - x^{-1}dx + \dots \in \Sw'(\R) \cdot dx =: \Omega^1 \Sw(\R) ,
\]
regarded as a tempered current. Then we have the \bigskip

{\bf Proposition}
\[
d\li_1(x) \; = \; \w(x) \cdot x^{-1} dx + L(x) dx \in \Omega^1_\Q \Q[[(x]] \;,
\]
\ie
\[
\frac{1}{1 + e^{-x}} + \frac{1}{1 - e^{-x}} =  \frac{2}{1 - e^{-2x}} .
\]
$\Box$ \bigskip


{\bf 2.2 Some Gel'fand - Shilov meromorphic functions} \bigskip

The holomorphic (in $s$) tempered (in $x$) distribution-valued function [GF I Ch 5.3]
\[
\gamma^s_+ (x)  := \frac{x_+^{s-1}}{\Gamma(s)}
\]
represents fractional differentiation $\partial^{-s}_x$ as convolution on smooth fast-decreasing functions, defining a group action $\gamma^s_+ * \gamma^t_+ = \gamma^{s+t}_+$ on generalised functions such as 
\[
\li_0(x) = (e^{-x}-1)^{-1} = - x^{-1} \sum_{n \geq 0} B_n \frac{x^n}{n!} \;.
\]
In \cite{4}(Prop 3) the tempered distribution  
\[
\li_s :=  \partial^{-s} \li_0 = \gamma^s_+ * \li_0
\]
(extending the function on the negative real axis defined by $\sum_{n \geq 1} e^{nx}n^{-s}$) is shown to satisfy the congruence
\[
\li_s \; \equiv  \;  - \pi \cot \pi s \cdot \gamma_+^s
\]
modulo smooth functions on the projective line. 

[The complicated term in our Proposition 3 simplifies, using [GF Ch I \S 4.4]
\[
(x \pm i0)^{s-1} = x_+^{s-1} \mp i e^{\pm i \pi s} x_-^{s-1} 
\]
to
\[
- \half [e^{-i \pi s} (x_+^{s-1} +  i e^{+ i \pi s}  x_-^{s-1})  \;   +  e^{+i \pi s}(x_+^{s-1} - i e^{-i \pi s}x_-^{s-1})]   \; = \;  -  \cos \pi s \cdot x_+^{s-1}.
\]
$\Box$ ]

This suggests defining a meromorphic regularization with smooth Laurent/Taylor coefficients [GF I \S 3.7] 
\[
\sli_s(x) := \li_s + \pi \cot \pi s \cdot \gamma_+^s \in C^\infty(\R)\{s\}
\]
with for example 
\[
 \sim  (1-s)^{-1} + \log \; |\frac{x}{e^x - 1}|  (1-s)^0  + \dots  \in \Q(((1-s)))
 \]
at $x=0$. \bigskip

 $\bullet$ I propose to regard $\R$ as a not-quite-probability space with $d\li_1$ as a not-quite distribution-valued measure, with  the remarkable (to me) property that its divided moments 
\[
\gamma^s_+ * d\li_1 := \gamma^s_+ *\li_0 dx = \li_s dx =  d\li_{s+1} \in \Omega^1\Sw(\R) \cong \Sw'(\R) dx
\]
are tempered. I contend below that it is reasonable in some physical models to think of the parameter $x$ as entropy.\bigskip

{\bf Proposition}  (\cite{4}(Prop 3 Cor 2)) {\it When $s=n$ is a positive integer, the $n$th divided moment of $x$ with respect to $d\li_1$ is congruent to
\[ 
\gamma^s_+ * d\li_1 \equiv -\gamma^{n-1}(x) \log |x| dx
\]
mod smooth functions. When $n=0$ we have $-x^{-1}dx$ as above. } \bigskip

For example, the tempered one-form $\gamma^s_+ * d\li_1$ is not integrable at $s=1$, where
\[
\int_\R \frac{x_+^{s-1}}{\Gamma(s)} \frac{1}{e^x - 1} \; dx  
\]
diverges. \bigskip

{\bf 2.3 Claim} 

In the language of Federer-de Rham distributions \cite{5,6,7},\cite{9}(\S 3.1) our proposition asserts that $d\li_1$ on the real projective line is congruent mod smooth measures to a logarithmic current, \ie the unit $\w$ times the Mellin measure $x^{-1}dx := d\mu(x)$ on the multiplicative group. It is a section of a line bundle, integral in the sense of both algebra [defined over $\Z$]  and arithmetic geometric measure theory. Here we invoke the fact that $\w(x) \in \Zt[[x]]$ is 2-adically integral and thus integrally flat in the sense of \cite{5}$\otimes \Zt$.

[Currents on the real line define a resolution \cite{13}  of $\R$ by acyclic sheaves \cite{6}(Theorem A.2), \cite{12}. The fundamental cohomology class $t^{-1}dt$ of the real projective line is Poincar\'e dual to a locally flat integral current.]\bigskip

{\bf Proposition} {\it The current $d\li_1$ satisfies the spark equation}
\[
da = \phi - R \in \Omega^1\Sw(\R)
\]

[As in \cite{5}(\S 0.1), with
\[
a = \li_1,  R = \w \cdot d\mu .
\] 

The equivalence class $[d\li_1] \in \hl^0(\bP_1(\R); \Zt)$ in a topological group of $\Zt$-adic Harvey-Lawson sparks is a generator.]\bigskip
 
{\bf 2.4} Planck's asymptotic expression \begin{footnote}{$h$ is Planck's constant, $k$ is Boltzmann's, and $c=1$; $\kappa = kh^{-1} \sim 2 \times 10^{11}$ Herz per Kelvin is a `coupling constant' used to express measurements of action in terms of units of heat \cite{26}.

The relation of heat to action is a question of the nature of matter, of measurement and physics, but it concerns ratios of ratios, arguably of the fourth order.}\end{footnote}
\[
\bro := \rho(\nu,T) \cdot d\nu  \sim  2h \frac{\nu^3} {e^{h\nu/kT} - 1} \cdot d\nu \in \Omega^1_\R \R[[\nu]]
\]
for the intensity of the black body radiation at frequency $\nu$ of a photon gas at temperature $T$ simplifies when written in terms  of the dimensionless ratio
\[
\ve = \frac{h\nu}{kT} = \kappa^{-1} \nu/T,  \ie \; \nu = \kappa T \ve
\]
to a version

\[
T^{-4} \bro  \sim 12h \kappa^4 \cdot \gamma^4_+*  d\li_1 \in \Omega^1_\R \R[[\ve]]
\]

of the Stefan-Boltzmann law in statistical mechanics: the radiance of a star scales as the fourth power of its temperature. Reinterpreted, this is a definition of temperature in terms of the fourth moment of entropy
with respect to the spark intensity $d\li_1$. 

{\bf 2.5} [digressive rant]: The resulting appearance of $\zeta(4)$ in the expression for the luminosity is particularly intriguing because Planck's photon has no underlying physical model: it is `matter' only in some rhetorical sense. This seems to display something very basic about the relation of heat to action in statistical mechanics.]\bigskip

 As for the second moment (or standard deviation or kinetic energy), the integrand of
\[
\gamma^2_+ * d\li_1 =  x \frac{1}{e^x - 1} dx 
\]
is the Hirzebruch multiplicative series \cite{15}
\[
\frac{x}{\exp_{{\mathbb G}_m}(x)}  = e^{-x/2} \frac{x/2}{\sinh x/2} =: e^{-x/2} \hat{A}(x)  
\]
for the arithmetic (Todd) genus of algebraic geometry, generalized by the Dirac operator to Spin manifolds \cite{27}(\S 2). \bigskip

{\bf Corollary} to Wu's lemma : {\it The characteristic classes defined for SU-bundles by the multiplicative series
\[
\hat{A}(2x) := (\w(x)\w(-x))^{1/2} = (\frac{2x}{e^{2x}-1} \cdot \frac{-2x}{e^{-2x}-1})^{1/2} \in \Zt [[x]] 
\]
reduce modulo two to the reciprocal
\[
(\nu^\spin_x)^{-1} \in \F_2[[x]] 
\]
of Hopkins and Singer's integral Wu classes for Spin bundles }\cite{9}(Appendix E.2).\bigskip

[These classes are defined by the multiplicative sequence  $(f(x)f(-x))^{1/2}$ with $f(x)^{-1} = \sum_{k \geq 0} x^{2^k-1}$. ]\bigskip

This is a developing topic (see \cite{22}(Part II)), which I hope to continue pursuing. These observations perhaps suggest a Schur-Weyl duality \cite{14, 21, 23} between the commutative monoid of points in a space and Boltzmann/Feynman/Witten statistical mechanics of the corresponding configuration spaces \cite{25} \dots \bigskip

{\bf \S III Hypothesis/Conjecture} [For Taira Honda, Yuri Manin, John McKay] \bigskip

Over $\R$ the Borel/Tate correspondence $ t_\T K((c)) \sim t_\T K((h))$ is the consequence of a congruence modulo smooth functions, defined by a product with the Wu \cite{28} spark \ $[d\li_1] \in \hl^0(\bP_1(\R); \Zt)$ \dots \bigskip

\bibliographystyle{amsplain}

\end{document}